\documentclass{article}

\usepackage[english]{babel}

\usepackage[letterpaper,top=2cm,bottom=2cm,left=3cm,right=3cm,marginparwidth=1.75cm]{geometry}

\usepackage{amssymb}
\usepackage{latexsym}
\usepackage{dsfont}
\usepackage{amsmath}
\usepackage{graphicx}
\usepackage[colorlinks=true, allcolors=blue]{hyperref}

\allowdisplaybreaks[4]               

\newtheorem{defn}{Definition}[section]
\newtheorem{thm}{Theorem}[section]

\newtheorem{lem}[thm]{Lemma}
\newtheorem{rem}{Remark}

\numberwithin{equation}{section}

\newcommand{\rmd}{\mathrm{d}}

\newcommand{\Var}{{\mathbb{V}}{\rm ar}}

\newcommand{\Enum}{\mathbb{E}}

\newcommand{\Pnum}{\mathbb{P}}

\newcommand{\Rnum}{\mathbb{R}}

\newcommand{\Nnum}{\mathbb{N}}

\newcommand{\lam}{\lambda}
\newcommand{\Lam}{\Lambda}
\newcommand{\vep}{\varepsilon}
\newcommand{\qed}{\hfill $\Box$}

\newcommand{\wt}{\widetilde}

\newcommand{\Dto}{\stackrel{{\rm d}}{\longrightarrow}}

\newcommand{\IID}{\stackrel{{\rm i.i.d.}}{\sim}}

\newcommand{\Proof}{ \noindent {\bf Proof.}\ }

\title{CLT, MDP and LDP for Range-Renewals of I.I.D.Samplings from an Infinite Discrete Distribution\footnote{Emails: chenxinx@sjtu.edu.cn; jsxie@fudan.edu.cn; zhaomz@zju.edu.cn;}}
\author{Xin-Xing Chen$^{1}$ \quad Jian-Sheng Xie$^{2}$ \quad Min-Zhi Zhao$^{3}$\\
{\footnotesize 1. School of Mathematical Sciences, Shanghai Jiaotong University, Shanghai, China}\\
{\footnotesize 2. School of Mathematical Sciences, Fudan University, Shanghai 200433, China}\\
{\footnotesize 3. School of Mathematical Sciences, Zhejiang University, Hangzhou, China}\\} 
\date{}

\begin{document}
\maketitle

\begin{abstract}
Let $R_n$ be the number of distinct values of the $n$ simple samples from an infinite discrete distribution. In 1960
Bahadure proved $\displaystyle \lim_{n\to \infty} \frac{R_n}{\Enum R_n}=1$ in probability; Chen et al. proved the limit 
in the sense of almost sure convergence, along with other results. In this note we present results of CLT, MDP and 
LDP for $R_n$ under mild conditions.

\end{abstract}

\section{Introduction}
Let $\pi$ be a probability measure on the set of natural numbers $\Nnum$ with $\pi_{_i} \ge \pi_{i+1}>0$ for all 
$i \geq 1$. Let $\xi:=\{ \xi_n \}_{n=1}^\infty$ be a sequence of independent random variables with common law 
$\pi$. Write
$$
R_n:=\#\{\xi_k: 1\le k\le n\}
$$
for the number of distinct values among the first $n$ elements of the process $\xi$. As in \cite{CXY} $\{R_n: \, 
n\ge 1\}$ is called the range-renewal process of $\xi$ . In 1960, Bahadur \cite{Bahadur} proved $\displaystyle 
\lim_{n\to \infty} \frac{R_n}{\Enum R_n}=1$ in probability. Recently, along with other results, Chen, Xie and Ying 
\cite{CXY} strengthened this limit to the sense of almost sure convergence, which could be regarded somehow as a 
strong law of large numbers (abbr. \textbf{SLLN}); Based on the work \cite{CXY}, Wu and Xie \cite{WX} discovered a delicate range-renewal structure 
in continued fractions. In this short note, we will go further and present results of central limit theorem (abbr. 
\textbf{CLT}), moderate deviation principle (abbr. \textbf{MDP}) and large deviation principle (abbr. \textbf{LDP}) 
for the quantity $R_n$ of i.i.d. samples $\{\xi_k\}_{k=1}^\infty$ from an infinite discrete distribution under mild 
conditions.

For simplicity of presentation, we set
\begin{equation}\label{e:defn_mu_sigma}
\mu(t) := \sum_{i=1}^\infty [1-e^{-t\pi_{_i} }] \text{ and } \sigma_{t} := \sqrt{\mu(2t)-\mu(t)}, \quad t\ge 0.
\end{equation}
As $n \to \infty$, we always have $\Enum R_n=\mu(n) \cdot [1+o(1)]$; under the mild conition $\displaystyle \varliminf_{n 
\to \infty} \Bigl( n^p \cdot \sigma_{n}^2 \Bigr)=\infty$ for some $p>0$ we also have $\Var(R_n)=\sigma_{n}^2 \cdot [1+o(1)]$; 
See Lemma \ref{lem: estimate4exp&var} for details. In general we have the following CLT result for $R_n$.
\begin{thm}\label{thm: CLT}
If $\sigma_{n} \to \infty$, then $\displaystyle \frac{R_n -\mu (n)}{\sigma_{n}} \Dto N (0,1)$.
\end{thm}

As the general LDP theory indicates \cite{D&Z}, we need more assumptions on the common law $\pi$ for the 
goal of establishing the MDP and LDP results for $R_n$. For the convenience of readers, we restate here the 
definitions of $\gamma$-regular functions and  $\gamma$-regular distributions introduced in \cite{CXY}.
\begin{defn}
Let $\zeta: [0,\infty) \to  [0,\infty)$ be a strictly increasing function. Let $\gamma\in [0,1]$. We
say that $\zeta$ is $\gamma$-regular if for any $\lam  >0$
\begin{equation}\label{eq: C2'-1}
\lim_{x \to \infty} \frac{\zeta (\lam x)}{\zeta (x)}=\lam^\gamma.
\end{equation}
\end{defn}

\begin{defn}\label{defn_pi_regular}
Let $\gamma\in (0,1)$. A distribution $\pi$ on $\Nnum$ is said to be $\gamma$-regular if there is a
$\gamma$-regular function $\zeta$ satisfying
\begin{equation}\label{eq:def-phi}
\lim_{n \to \infty} \pi_n \cdot \zeta^{-1}(n)  =1,
\end{equation}
where $\zeta^{-1}$ is the inverse function of $\zeta$.
\end{defn}

We refer the readers to \cite{CXY} for the definition of $1$-regular distributions. For any $\gamma$-regular
distribution $\pi$ on $\Nnum$ with $\gamma \in (0,1]$, we will call $\gamma (\pi) :=\gamma$ the regular
index of $\pi$. The above concepts are closely related to the concept of regular varying functions introduced
by Karamata \cite{Karamata}; For related properties of such functions, one is referred to \cite{BGT} .

For any set $\Gamma \subset \Rnum$, let $\Gamma^o$ and  $\bar{\Gamma}$  denote the interior and the closure of $\Gamma$ 
respectively. Our MDP and LDP results for $R_n$ are presented respectively as the following two theorems.
\begin{thm}\label{thm: MDP}
Let $b (n)$ be positive such that $b(n) \to \infty$ and $\frac{\mu (n)}{b (n)} \to \infty$ as $n \to \infty$. 
Let 
$$
Z_n:=\frac{R_n-\mu(n)}{\sqrt{\mu(n)b(n)}}.
$$ 
If $\pi$ is regular with index $\gamma \in (0,1)$, then $\{Z_n\}_{n=1}^\infty$ satisfies an LDP with rate 
function $I_\gamma (x)=\frac{x^2}{2 (2^{^\gamma}-1)}$. More precisely, for every measurable set $\Gamma 
\subset \Rnum$
\begin{eqnarray*}
-\frac{1}{2 (2^{^\gamma}-1)} \inf_{x \in \Gamma^o} x^2 \leq& \varliminf \limits_{n \to \infty} \frac{1}{b (n)} 
\log \Pnum (Z_n \in \Gamma)&\\
\leq& \varlimsup \limits_{n \to \infty} \frac{1}{b (n)} \log \Pnum (Z_n \in \Gamma)&\leq -\frac{1}{2 
(2^{^\gamma}-1)} \inf_{x \in \bar{\Gamma}} x^2.
\end{eqnarray*}
\end{thm}

\begin{thm}\label{thm: LDP}
Let $Y_n=\frac{R_n-\mu(n)}{\mu(n)}$. If $\pi$ is regular with index $\gamma\in (0, 1)$, then
$\{Y_n\}_{n=1}^\infty$ satisfies the following LDP: for every measurable set $\Gamma \subset \Rnum$
\begin{eqnarray*}
-\inf_{x \in \Gamma^o} \Lam^*_{\gamma} (x) \leq& \varliminf \limits_{n \to \infty} \frac{1}{\mu (n)} \log
\Pnum (Y_n \in \Gamma)&\\
\leq& \varlimsup \limits_{t \to \infty} \frac{1}{\mu (n)} \log \Pnum (Y_n \in \Gamma)&\leq -\inf_{x \in
\bar{\Gamma}} \Lam^*_{\gamma} (x),
\end{eqnarray*}
where $\Lam^*_{\gamma}(x)=\sup\limits_{\lam \in  R}\{\lam   x-\Lam_{\gamma}(\lam  )\}$ is the Fenchel-Legendre
transform of a smooth function $\Lam_{\gamma}$ given by
\begin{equation}\label{eqno: 0-def4Lam}
\Lam_{\gamma}(\lam  )=-\lam  +\frac{\gamma}{\Gamma(1-\gamma)}\int_0^\infty \log \bigl[ 1+(e^\lam  -1)(1-e^{-s})
\bigr] \cdot \frac{\rmd s}{s^{1+\gamma}}.
\end{equation}
The function $\Lam_{\gamma}$ has an explicit expansion for $\lam<\log 2$
\begin{equation}\label{eqno: 1-def4Lam}
\Lam_{\gamma} (\lam)=-\lam+\sum_{n=1}^\infty (1-e^{\lam})^n \cdot \sum_{k=1}^n \frac{(-1)^k}{k^{1-\gamma}}
C_{n-1}^{k-1}.
\end{equation}
\end{thm}

\begin{rem}
When $\pi$ is 1-regular, one may have  MDP and LDP for $R_n$ with respectively $I_1 (x)=\frac{x^2}{2}$ and 
$\Lam _1(\lam  )=e^\lam  -1-\lam$. We omite the corresponding proof in this note since it is much similar to 
the case $\gamma \in (0,1)$.
\end{rem}

The rest of the paper is organized as follows. In Section \ref{sec: Coupling}, we coupling the process $\xi$ with a standard
poisson process $\{N_t: \, t\ge 0\}$, so that $R^*_t:=R_{N_t}$ has a simple decomposition as an independent sum. In Section 
\ref{sec: Estimations} we present some preliminary estimations. With the help of $R^*_t$ and related estimations, we present 
the proofs of Theorems \ref{thm: CLT}--\ref{thm: LDP} in the last two sections, Sections \ref{sec： proof4CLT} and 
\ref{sec： proof4MDP&LDP}.
\section{Coupling the Process $\xi$ with a Standard Poisson Process}\label{sec: Coupling}

Let $\{N_t: t \geq 0\}$ be a standard Poisson process independent of $\xi=\{\xi_i\}_{i=1}^\infty \IID \pi$ and set 
\begin{equation}
N^{(i)}_t :=\sum_{k=1}^{N_t} 1_{\{ \xi_k =i \}}, \quad i\ge 1, ~t\ge 0.
\end{equation}
By the general theory of Poisson processes, $N^{(i)} :=\{N^{(i)}_t: t \geq 1\}$ is also a time-homogeneous Poisson
process with intensity $\pi_{_i} $. Moreover, these $\{N^{(i)}\}_{i=1}^\infty$ are independent processes.

Define $R^*_t:=R_{N_t}=\#\{\xi_k: 1 \leq k \leq N_t\}$. Then $R^*_t$ can be rewritten as
\begin{equation}\label{eq: dec4R}
R^*_t=\sum_{i=1}^\infty 1_{\{N^{(i)}_t \geq 1\}}.
\end{equation}
The right hand of \eqref{eq: dec4R} is an independent sum due to the independence of $\{N^{(i)}\}_{i=1}^\infty$. 
In view of \eqref{eq: dec4R} our problems under investigation will become simpler when turning to $R^*_t$ instead 
of $R_n$. Under some mild condtions we shall prove that $R^*_n-R_n$ is negligible in the scale of $\sigma_{n}$.

\section{Preliminary Estimates}\label{sec: Estimations}
In this part we will show how $R_n^*$  is close to $R_n$. But before doing this, we need some estimates on the 
first and second order derivatives of $\mu(t)$, denoted by $\dot{\mu}$ and  $\ddot{\mu}$ respectively. Clearly
\begin{equation}\label{eq: def4derivatives}
\dot{\mu} (t)=\sum_{i=1}^\infty \pi_{_i}  e^{-t \pi_{_i} }, \; \quad  \ddot{\mu} (t)=-\sum_{i=1}^\infty \pi_{_i} ^2 e^{-t \pi_{_i} }.
\end{equation}
It's obvious that $0<\dot{\mu} (t)<1$ and $-1<\ddot{\mu}(t)<0$ for all $t>0$.

\begin{lem}\label{lem:preparation}
{\rm(1)} For each $t>0$, we have
\begin{equation}\label{eq: error-bound-2}
\frac{\sigma_{t}^2}{t} \leq \dot{\mu} (t) \leq \frac{\mu (t)}{t}, \quad
-\ddot{\mu} (t) \leq \frac{2\mu (t)}{t^2}.
\end{equation}
In addition,
\begin{equation}\label{eq: bound-4-1-derive}
\dot{\mu} (s) \leq \bigl[ \dot{\mu} (t) \bigr]^{s/t}, \quad \hbox{for all } 0<s<t.
\end{equation}

{\rm (2)} $\sigma_{t}^2\le \mu(t) \le t$ for all $t>0$. Furthermore, $\mu(t) = o(t)$ and $\sigma_{t}^2=o(t)$,  
as $t \to \infty$.

{\rm (3)} For any $0<\vep<1$ and $t>0$,
\begin{equation}\label{eq: var-vs-dotmu}
\dot{\mu} (t) \le \Bigl(\frac{1}{\vep} +2\Bigr) \cdot \Bigl( \frac{\sigma_{t}^2}{t} \Bigr)^{1-\vep}, \quad -\ddot{\mu} (t)
\le \frac{1}{e\varepsilon} \cdot \frac{\bigl[\dot{\mu} (t) \bigr]^{1-\vep}}{t}.
\end{equation}
\end{lem}
\noindent\textbf{Proof.\;}
(1) Since  $\ddot{\mu} (t)<0$ for all $t> 0$, $\dot{\mu} (t)$ is strictly decreasing in $t$. And we have 
$$ 
\dot{\mu} (t)\ge \frac{\mu(2t)-\mu(t)}{t}=\frac{\sigma_{t}^2}{t}.
$$
Noting the simple fact $\frac{x^k}{k!} \cdot e^{-x}\le 1-e^{-x}$ for all   $k\ge 1$ and all $x>0$, we have
$$
\sum_{i=1}^\infty \pi_{_i} ^{k} \cdot e^{-t\pi_{_i} }\le \sum_{i=1}^\infty \frac{k!}{t^k} \cdot (1-e^{-t\pi_{_i} }) 
=k! \cdot \frac{\mu(t)}{t^k}.
$$
Combining \eqref{eq: def4derivatives} with the above estimation, we easily prove the rest two inequalities in 
\eqref{eq: error-bound-2}.

In view of $\xi_1 \sim \pi$ and $\dot{\mu} (t)=\Enum [\exp(-t  \pi_{_{\xi_1}})]$, \eqref{eq: bound-4-1-derive} follows 
immediately from Jensen's inequality.

(2) The inequality $\sigma_{t}^2 \leq \mu(t)$ is already contained in \eqref{eq: error-bound-2}. Applying $1-e^{-x} 
\leq x,  \forall x \geq 0$ in Eq. \eqref{e:defn_mu_sigma}, we have $\displaystyle \mu (t) \leq \sum_{i=1}^\infty t \pi_{_i} =t$. 
In view of Lebesgue's dominated convergence theorem 
$$
\lim\limits_{t\to \infty}\frac{\mu(t)}{t}=\lim_{t\to \infty} \sum\limits_{i=1}^\infty \frac{1-e^{-t\pi_{_i} }}{t\pi_{_i} } \cdot \pi_{_i} = \sum_{i=1}^\infty
\lim_{t\to \infty} \frac{1-e^{-t\pi_{_i} }}{t \pi_{_i} } \cdot \pi_{_i} =0,
$$
which implies $\mu (t) = o(t)$ and hence $\sigma_{t}^2=o(t)$.

(3)  Fix  $0<\vep<1$. Noting  $1-e^{-x}\ge \frac{x}{2}$ for all $x\in  [0,1]$ and $0<\frac{\sigma_{t}^2}{t} \leq 1$, we have
\begin{equation}\label{e:pi1up}
\sum_{i: t \pi_{_i}  \le 1} \pi_{_i}  e^{-t \pi_{_i} } \le  \frac{2}{t} \sum_{i: t \pi_{_i}  \le 1} e^{-t \pi_{_i} } \cdot  (1-e^{-t \pi_{_i} })
\le 2\cdot \frac{\sigma_{t}^2}{t}\le 2\cdot\Bigl( \frac{\sigma_{t}^2}{t} \Bigr)^{1-\vep}.
\end{equation}
Let $p:=1/\vep-1>0$. By $\displaystyle \max_{x\ge 0} \bigl( x^p e^{-x} \bigr) =p^p e^{-p}$, we can write
\begin{equation}
\sum_{i=1}^\infty \pi_{_i} ^{p+1} e^{-t \pi_{_i} } \le \frac{p^p e^{-p}}{t^p} \sum_{i=1}^\infty  \pi_{_i}  =
\Bigl( \frac{p}{et} \Bigr)^p.
\end{equation}
In view of H\"{o}lder's inequality we obtain
\begin{align*}
\sum_{i:t \pi_{_i} >1} \pi_{_i}  e^{-t \pi_{_i} } \le& \Bigl[ \sum_{i:t \pi_{_i} >1} \pi_{_i} ^{p+1} e^{-t \pi_{_i} } \Bigr]^{1/(p+1)} \cdot \Bigl[
\sum_{i:t \pi_{_i} >1} e^{-t \pi_{_i} }\Bigr]^{p/(p+1)}\\
\leq& \Bigl( \frac{p}{et} \Bigr)^{p/(p+1)} \cdot \Bigl[ \sum_{i:t \pi_{_i} >1} e^{-t \pi_{_i} } \cdot \frac{1-e^{-t \pi_{_i} }}{1-e^{-1}}
\Bigr]^{p/(p+1)}\\
= & \Bigl(\frac{p}{e-1}\Bigr)^{p/(p+1)}  \cdot \Bigl[\frac{1}{t} \sum_{i:t \pi_{_i} >1} e^{-t \pi_{_i} }(1-e^{-t \pi_{_i} }) \Bigr]^{p/(p+1)}\\
\leq & \Bigl(\frac{p}{e-1}\Bigr)^{p/(p+1)} \cdot  \Bigl(\frac{\sigma_{t}^2}{t} \Bigr)^{p/(p+1)}= \Bigl[\frac{1-\vep}{(e-1)\vep} 
\Bigr]^{1-\vep}\cdot \Bigl( \frac{\sigma_{t}^2}{t} \Bigr)^{1-\vep}\\
\leq & \frac{1}{\vep}\cdot \Bigl( \frac{\sigma_{t}^2}{t} \Bigr)^{1-\vep}.
\end{align*}
Combining the above estimate with \eqref{e:pi1up}, we complete the proof of the first inequality in \eqref{eq: var-vs-dotmu}.

Similarily we can apply H\"{o}lder's inequality with $p=\frac 1{\varepsilon}$ and $q=\frac p{p-1}$, yielding
\begin{eqnarray*}
\sum_{i=1}^\infty \pi_{_i} ^2 e^{-t \pi_{_i} } 
&\le& \Bigl[ \sum_{i=1}^\infty \pi_{_i} ^{p+1} 
e^{-t \pi_{_i} } \Bigr]^{1/p} \cdot \Bigl[ \sum_{i=1}^\infty \pi_{_i}  e^{-t \pi_{_i} } \Bigr]^{1/q}\\
&\leq&  \frac {p}{e t}\Bigl[ \sum_{i=1}^\infty \pi_{_i}  e^{-t \pi_{_i} } \Bigr]^{1/q}=\frac {1}{ e \vep} \cdot \frac{[\dot{\mu}(t)]^{1-\vep}}{t}.
\end{eqnarray*}
This proves the second inequality in \eqref{eq: var-vs-dotmu}.
\qed

\begin{lem}\label{lem: Difference2CLT}
We always have $\frac{R_{n}^*-R_n}{\sigma_{n}} \to 0$ in probability as $n \to \infty$.
\end{lem}
\noindent\textbf{Proof.\;}
 By Lemma \ref{lem:preparation}, $\frac{\mu(n)}{n} \to 0$. Let $w_n=[\frac{n^4}{\mu (n)}]^{1/6}$,  then 
$$
\frac{w_n}{n}= \Bigl[ \frac{1}{n^2 \cdot \mu (n)} \Bigr]^{1/6} \to 0 \hbox{ and } 
\frac{n }{w_n^2} =\bigl[ \frac{\mu (n)}{n} \bigr]^{1/3} \to 0.
$$
Hence there exists some $n_0\ge 1$ with $w_n/n<1/6$ for all $n \ge n_0$. Noting $0<\dot{\mu}(t)<1, 
\ddot{\mu}(t)<0$ and inequalities \eqref{eq: error-bound-2}--\eqref{eq: bound-4-1-derive}, we can write for all 
$n \ge n_0$
\begin{eqnarray*}
\nonumber &&\frac{\mu (n+w_n)-\mu (n-w_n)}{\sigma_{n}}\leq \frac{2w_n}{\sigma_{n}}  \cdot \dot{\mu}(n-w_n) 
\leq \frac{2w_n}{\sigma_{n}} \cdot \dot{\mu} (n)^{1-w_n/n}\\
\nonumber  &\leq&  \frac{2w_n}{\sigma_{n}} \cdot \dot{\mu} (n)^{5/6} 
\leq \frac{2w_n}{\sigma_{n}} \cdot  \Bigl( \frac{\sigma_{n}^2}{n} \Bigr)^{5/6} \leq 2 \Bigl[ \frac{n^4}{\mu (n)} 
\Bigr]^{1/6} \cdot \frac{\sigma_{n}^{4/6}}{n^{5/6}}=2 \cdot \Bigl[\frac{\sigma_{n}^4}{n \mu(n)}\Bigr]^{1/6}.
\end{eqnarray*}
Noting $\sigma_{n}^2 \leq \mu (n)$ and $\frac{\mu (n)}{n} \to 0$ as $n \to \infty$, we obtain
\begin{equation}\label{e:+w-wupperbound}
\frac{\mu (n+w_n)-\mu (n-w_n)}{\sigma_{n}} \leq 2 \cdot \Bigl[\frac{\mu(n)}{n}\Bigr]^{1/6} \to 0.
\end{equation}

Now fix $\vep >0$. In view of $R^*_t=R_{N_t}$, we clearly have $R_{n-w_n}^*\le R_n^*\le R_{n+w_n}^*$ and
$$
\{N_{n-w_n}\le n\le N_{n+w_n}\}\subseteq \{R_{n-w_n}^*\le R_n \le R_{n+w_n}^*\}
$$
If $|R_n^*-R_n|\ge \vep  \sigma_{n}$ and $N_{n-w_n}\le n\le N_{n+w_n}$, then $
R_{n+w_n}^*-R_{n-w_n}^*\ge \vep  \sigma_{n}$. It follows immediately
\begin{equation}\label{e:minuR*Rscale}
\Pnum(|R_n^*-R_n|\ge \vep  \sigma_{n})\le \Pnum( R_{n+w_n}^*-R_{n-w_n}^*\ge \vep  \sigma_{n})+\Pnum(N_{n-w_n}>n)+\Pnum(N_{n+w_n}<n).
\end{equation}
\eqref{e:+w-wupperbound} implies $\Pnum( R_{n+w_n}^*-R_{n-w_n}^*\ge \vep  \sigma_{n})\le \frac{\mu_{n+w_n}-\mu_{n-w_n}}{\vep  
\sigma_{n}} \to 0$. Noting $\Pnum(|N_s-s|\ge x) \le \frac{s}{x^2}$ for  all $s,x>0$, we also have
\begin{eqnarray*}
&& \Pnum(N_{n-w_n}>n)+\Pnum(N_{n+w_n}<n)\\
&=& \Pnum(N_{n-w_n}-(n-w_n)>w_n)+\Pnum((n+w_n)-N_{n+w_n}>w_n)\\
&\le& \frac{n-w_n}{w_n^2}+\frac{n+w_n}{w_n^2}= \frac{2n}{w_n^2} \to 0.
\end{eqnarray*}
Hence $\displaystyle \lim_{n\to \infty} \Pnum(|R_n^*-R_n|\ge \vep  \sigma_{n})=0$ for each $\vep >0$, 
proving the lemma.
\qed

Next, we calculate the moments of $R_n^*$ and $R_n$, and compare their differences. Clearly
\begin{equation}
\Enum R_n^* = \sum_{i=1}^\infty (1-e^{n\pi_{_i} })=\mu(n), \quad \Var(R_n^*)=\mu(2n)-\mu(n).
\end{equation}
By setting for each $n\ge 1$
\begin{equation}
\Delta(n) :=2 \sum_{(i,j): j>i\ge 1} \bigl[ (1-\pi_{_i} )^n (1-\pi_{_j})^n-(1-\pi_{_i} -\pi_{_j})^n \bigr],
\end{equation}
we also have
\begin{equation}
\Enum R_n = \sum_{i=1}^\infty \bigl[ 1-(1-\pi_{_i} )^n \bigr], \quad \Var(R_n) =  \Enum R_{2n} -\Enum R_n-\Delta(n).
\end{equation}

\begin{lem}\label{lem: estimate4exp&var}
{\rm (1)} For any $n \geq 1$,
\begin{equation}
\label{eq: error-bound-1} \mu(n) \le \Enum  R_n   \le  (1+\frac{e}{n}) \cdot  \mu(n).
\end{equation}

{\rm (2)} If $\displaystyle \lim_{n \to \infty} [n^p \cdot \sigma_{n}^2]=\infty$ for some $p>0$, then as $n\to \infty$,
\begin{eqnarray}\label{e:VarR*Rcom}
\Var(R_n)=\sigma_{n}^2 \cdot [1+o(1)].
\end{eqnarray}

{\rm (3)} If $\pi$ is regular with index $\gamma\in (0,1)$ and $\zeta$ being the regular function defined in \eqref{eq:def-phi}, then
\begin{equation}\label{e:mu estimate}
\mu(t)=\Gamma(1-\gamma)\zeta(t) \cdot [1+o(1)], \quad \sigma_{t}^2=(2^{^\gamma}-1) \Gamma(1-\gamma)\zeta(t) \cdot [1+o(1)],\quad t\to \infty.
\end{equation}
\end{lem}
{\bf Proof.} (1) Noting  $\Enum R_n-\mu(n) = \sum\limits_{i=1}^\infty [e^{-n \pi_{_i} } -(1-\pi_{_i} )^n]$ and the fact $e^{-x}\ge 1-x$ for all $x\ge 0$, 
we have the lower bound of \eqref{eq: error-bound-1}. Since $x^n-y^n\le n(x-y) x^{n-1}$ and $e^{-x}\le 1-x+\frac{x^2}{2}$ for all $x\ge  y\ge 0$ 
and all $n\ge 1$, we also have
\begin{equation}\label{e:minusexpRR*}
\Enum R_n-\mu(n) \le n \sum_{i=1}^\infty  (e^{-\pi_{_i} }-1+\pi_{_i} ) e^{-(n-1) \pi_{_i} } \le  n \sum_{i=1}^\infty \frac{\pi_{_i} ^2}{2}  e^{-n \pi_{_i} }\cdot 
\max\limits_{i\ge 1} e^{\pi_{_i} }\le \frac{n e}{2} \cdot |\ddot{\mu}(n)|.
\end{equation}
By \eqref{eq: error-bound-2}, $|\ddot{\mu}(n)| \le \frac{2 \mu(n)}{n^2}$. This proves the upper bound of \eqref{eq: error-bound-1}.

(2) Exploiting the fact $x^n-y^n\le n(x-y) x^{n-1}$  for all $x\ge  y\ge 0$ again, we can easily get
\begin{eqnarray*}
0 \leq \Delta (n) &\leq& 2\sum_{(i,j):j>i\ge 1} n  \pi_{_i} \pi_{_j} (1-\pi_{_i} -\pi_{_j}+\pi_{_i} \pi_{_j})^{n-1}\\
&\leq& n \Bigl[ \sum_{i=1}^\infty \pi_{_i}  (1-\pi_{_i} )^{n-1} \Bigr]^2 \leq \frac{n}{(1-\pi_{_1})^2} \Bigl( \sum_{i=1}^\infty \pi_{_i}  e^{-n \pi_{_i} }
\Bigr)^2,
\end{eqnarray*}
where the last inequality holds since $1>\pi_{_1}\ge \pi_{_i} $ for all $i\ge 2$. Applying \eqref{eq: var-vs-dotmu} with $\vep=\frac{1}{4}$,  
we get  $\sum\limits_{i=1}^\infty \pi_{_i}  e^{-n \pi_{_i} }=\dot{\mu} (n)=O \Bigl( \Bigl[ \frac{\sigma_{n}^2}{n} \Bigr]^{\frac 34} \Bigr)$. Noting $\frac{\sigma^2}{n} 
\leq \frac{\mu(n)}{n}=o(1)$, it follows immediately
$$ 
\Delta (n) =O \Bigl( n \Bigl[ \frac{\sigma_{n}^2}{n} \Bigr]^{\frac 32}\Bigr)=   O\Bigl( \sigma_{n}^2\cdot \Bigl[ \frac{\mu(n)}{n} \Bigr]^{\frac 12} \Bigr)=o(\sigma_{n}^2).
$$
Now suppose  $\displaystyle \lim_{n\to \infty} [n^p \cdot \sigma_{n}^2]=\infty$ for some $p>0$. Applying Lemma \ref{lem:preparation} (3) with $\vep=1-\sqrt{\frac{p}{p+1}}$, 
we have $|\ddot{\mu} (n)| =O(\frac{\dot{\mu}(n)^{1-\vep}}{n}) =O(\frac{1}{n} \cdot \Bigl[\frac{\sigma_{n}^2}{n} \Bigr]^{(1-\vep)^2})= O \Bigl( \frac{1}{n} \cdot \Bigl[ 
\frac{\sigma_{n}^2}{n} \Bigr]^{\frac{p}{p+1}}\Bigr)$. This, together with \eqref{e:minusexpRR*}, yields
$$
|\Enum R_n-\mu (n)|  =  O(n |\ddot{\mu} (n)|) = O \Bigl( \Bigl[ \frac{\sigma_{n}^2}{n} \Bigr]^{\frac{p}{p+1}}\Bigr) 
= O \Bigl( \sigma^2_n \cdot \frac{1}{(n^p  \cdot \sigma_{n}^2)^{\frac 1{p+1}}} \Bigr) = o(\sigma_{n}^2).
$$
Since $|\ddot{\mu} (2n)|\le |\ddot{\mu} (n)|$, we also have  $|\Enum  R_{2n}-\mu (2n)|=O(n |\ddot{\mu} (2n)|)=O(n |\ddot{\mu} (n)|)=  o(\sigma_{n}^2)$.
Therefore 
\begin{eqnarray*}
\left| \Var(R_n)-\Var(R_n^*) \right| &\leq& |\Delta (n)|+|\Enum R_n -\mu (n)|+| \Enum R_{2n} -\mu (2n)| =o(\sigma_{n}^2).
\end{eqnarray*}
This proves \eqref{e:VarR*Rcom}.

 (3) By \cite[Lemma 3.3]{CXY}, we have $\Enum R_n=\Gamma(1-\gamma)\cdot \zeta(n)\cdot [1+o(1)]$, yielding \eqref{e:mu estimate} immediately
with the results proved above.
\qed

\begin{rem}
Let $\pi_{_i}:=2 \bigl[ (\frac{1}{2})^{i!} -(\frac{1}{2})^{(i+1)!} \bigr]$ for all $i \geq 1$. This defines a distribution $\pi$ on $\Nnum$. Let $t_n:=2^{n! \cdot 
\sqrt{n+1}}$. We have $\pi[n+1,\infty)=2 \cdot t_n^{-\sqrt{n+1}}$ and $(\frac{1}{2})^{i!} \leq \pi_{_i} \leq (\frac{1}{2})^{i!-1}$ for all $i \geq 1$. Clearly
\begin{eqnarray*}
\sigma^2_{{t_{_n}}} &=& \sum_{i=1}^\infty e^{-t_n \cdot \pi_{_i}} \cdot [1- e^{-t_n \cdot \pi_{_i}}]
\leq n e^{-t_n \cdot \pi_{_n}} +\sum_{i=n+1}^\infty t_n \cdot \pi_{_i}\\
&\leq& n \exp\{-\frac{1}{2} \cdot t_n^{1-1/\sqrt{n+1}}\}+2 \cdot t_n^{1-\sqrt{n+1}}.
\end{eqnarray*}
This implies  $\displaystyle \varliminf_{n \to \infty} [n^p \cdot \sigma_{n}^2]=0$ for all $p>0$. Therefore the assumption in Lemma \ref{lem: estimate4exp&var} {\rm (2)} makes sense.
\end{rem}

At the end of this section, we give an estimate needed in the proof of Theorem \ref{thm: LDP} in Section \ref{sec： proof4MDP&LDP}.
\begin{lem}\label{l:Lambda_0.9}
Let $\zeta$ be a regular function with index $\gamma\in (0,1)$ and inverse function $\zeta^{-1}=\varphi$. Then 
$$
\lim_{t\to \infty} \frac{1}{\zeta(t)}\sum_{i=1}^\infty \log\bigl[1+(1-e^{-\frac{t}{\varphi (i)}}) u \bigr]
=\gamma\int_0^\infty \log\bigl[1+(1-e^{-s})u\bigr] \cdot \frac{\rmd s}{ s^{1+ \gamma}}, \quad \forall u>-1.
$$
\end{lem}
\Proof
The case $u=0$ is trivial. Fix $u>-1$ and $u \neq 0$. Let $t\ge 1$ and  write ${t_*}=\zeta(t)$ for short. We define a (discrete) measure $m_t$ on $(0,\infty)$ as below: 
for each $i\ge 1$,
$$
m_t((\frac{i-1}{{t_*}}, \frac{i}{{t_*}}])=m_t(\{\frac{i}{{t_*}}  \})=  \frac{1-e^{-\frac{t}{\varphi(i)}}}{{t_*}}.
$$
Put $\alpha:=1/\gamma$ and define also a measure $m$ on  $(0,\infty)$ as $m(\rmd  x)= (1-e^{-1/x^\alpha}) \rmd  x$. Define a function $f: (-1,\infty) \to \Rnum$ as the 
following
$$
f(x):=\left\{
\begin{array}{rl}
\frac{\log (1+x)}{x}, & \hbox{ if } x>-1 \hbox{ and } x \neq 0,\\
1, & \hbox{ if } x =0.
\end{array}
\right.
$$
Then we define for each $x>0$ $h_t(x):=f((1-e^{-\frac{t}{\varphi(x \cdot {t_*})}})u) \cdot u$ and $h(x):=f((1-e^{-1/x^{\alpha}})u) \cdot u$.
Clearly we have $\int_0^\infty h(x)m(\rmd  x)= \gamma\int_0^\infty \log\bigl[1+(1-e^{-s})u\bigr] \cdot \frac{\rmd s}{s^{1+ \gamma}}$ and
$$
\int_0^\infty h_t(x)m_t(\rmd  x) = \sum_{i=1}^\infty  h_t(\frac{i}{{t_*}})m_t(\{\frac{i}{{t_*}}\})=\frac{1}{\zeta(t)}\sum_{i=1}^\infty \log\bigl[1+(1-e^{-\frac{t}{\varphi (i)} })u\bigr].
$$
Thus in order to prove our lemma,  it suffices to show that
\begin{equation}\label{eqi:equivalentcondition}
\lim_{t\to\infty}\int_0^\infty h_t(x)m_t(\rmd  x)=\int_0^\infty h(x)m(\rmd  x).
\end{equation}

We seperate the proof of \eqref{eqi:equivalentcondition} into the following parts.

\textbf{(A)} First we prove that $m, m_t, t\ge 1$ are finite measures satisfying 
\begin{equation}\label{eq: finite-meas-1}
\lim_{t\to \infty} m_t((0, \infty))=m((0,\infty)).
\end{equation}
In fact, by Lemma 3.3 of \cite{CXY}, $\displaystyle \lim_{t \to \infty} m_t((0,\infty))=\lim_{t \to \infty}\frac{1}{{t_*}} \sum_{i=1}^\infty \bigl[ 
1-e^{-\frac{t}{\varphi(i)}} \bigr]=\Gamma(1-\gamma)$. On the other hand, a direct calculation indicates
\begin{align*}
m((0,\infty))=&\int_0^\infty   (1-e^{-x^{-1/\gamma}}) \rmd  x= -\int_0^\infty  x   \rmd (1-e^{-x^{-1/\gamma}})\\
=&  -\int_0^\infty x e^{-x^{-1/\gamma}} \rmd x^{-1/\gamma}=\int_0^\infty y^{-\gamma} e^{-y} \rmd y=\Gamma(1-\gamma),
\end{align*}
proving \eqref{eq: finite-meas-1}.

\textbf{(B)} Then we prove that for each $x >0$
\begin{equation}\label{eq: finite-meas-2} 
\lim_{t\to \infty}h_t(x)=h(x) \hbox{ and } \lim_{t\to \infty} m_t((0, x])=m((0,x]).
\end{equation}
Since $\zeta$ is $\gamma$-regular, in view of \cite[Theorem 1.5.1]{BGT} the limit
\begin{equation}\label{e:gammaregular}
\lim_{t\to \infty}\frac{t}{\varphi(x \cdot {t_*})}=x^{-1/\gamma}
\end{equation}
holds uniformly on each compact subset of $(0,\infty)$. Then the first limit in \eqref{eq: finite-meas-2} holds. Now Fix $x\in (0,\infty)$. 
By (\ref{e:gammaregular}) and Lebesgue's dominated convergence theorem,
\begin{align*}
\lim_{t\to \infty}\int_0^x \bigl[ 1-e^{-\frac{t}{\varphi(y{t_*})}} \bigr] \rmd  y=\int_0^x \lim_{t\to \infty}
\bigl[ 1-e^{-\frac{t}{\varphi(y{t_*})}} \bigr] \rmd  y=\int_0^x (1-e^{-y^{-1/\gamma}}) \rmd  y=m((0,x]).
\end{align*}
By the increasing property of $\varphi$,
\begin{equation}\label{e:m_tconvergence}
\int_0^{ x{t_*}}\frac{1-e^{-\frac{t}{\varphi(s)}}}{{t_*}} \rmd  s\ge  \sum_{i=1}^{\lfloor x{t_*}\rfloor}\frac{1-e^{-\frac{t}{\varphi(i)}}}{{t_*}}
\ge  \int_1^{  x{t_*}-1 }\frac{1-e^{-\frac{t}{\varphi(s)}}}{{t_*}} \rmd  s.
\end{equation}
Since both the upper and lower bounds in \eqref{e:m_tconvergence} converge to $\lim\limits_{t\to \infty}\int_0^x (1-e^{-\frac{t}{\varphi(y{t_*})}}) \rmd  
y=m((0,x])$ as $t\to \infty$, we have 
$$
\lim_{t\to \infty} m_t((0, x])=\lim_{t\to \infty} \sum_{i=1}^{\lfloor x{t_*}\rfloor}\frac{1-e^{-\frac{t}{\varphi(i)}}}{{t_*}}    
=\lim_{t\to \infty}\int_0^x (1-e^{-\frac{t}{\varphi(y{t_*})}}) \rmd  y=m((0,x]),
$$
proving the second limit in \eqref{eq: finite-meas-2}.

\textbf{(C)} Now we prove $\sup\{|h_t(x)|: x\in (0,\infty), ~t\ge 1\}<\infty$ and $\sup\{|h(x)|: x\in (0,\infty)\}<\infty$.
In fact, noting $0<|\bigl[ 1-e^{-\frac{t}{\varphi(x t_*)}} \bigr] \cdot u| \leq |u|$ for all $t\ge 1$ and $x>0$, we have
$$
|h_t(x)|\le C_u \cdot |u|,
$$
where $C_u= \max\{|\frac{\log(1+y)}{y}|: 0<|y| \leq |u| \}<\infty$. Similarly, we have $|h(x)|\le C_u \cdot |u|$ for all $x>0$.

Finally based on the above facts in parts \textbf{(A)}--\textbf{(C)}, one easily proves \eqref{eqi:equivalentcondition}.
\qed

\section{Proof of Theorem \ref{thm: CLT}}\label{sec： proof4CLT}
To prove  Theorem \ref{thm: CLT}, we need Esseen's inequality for normal approximation; The
following version of Esseen's inequality is in fact a direct consequence of \cite[Theorem 7.1]{CGS}.
\begin{lem}\label{lem: CLT}
Let $\{X_n\}_{n=1}^\infty$ be independent  random variables with $\Enum X_n=0$ for each $n$. Suppose
$$
\sum_{n=1}^\infty \Enum [X_n^2] =1, \; \sum_{n=1}^\infty \Enum [|X_n|^3] <\infty,
$$
hence $S=\sum_{n=1}^\infty X_n$ converges almost surely. Let $F(x)$ be the cumulative distribution function of $S$ and let $\displaystyle \Phi (x):=\int_{-\infty}^x 
\frac{e^{-t^2/2}}{\sqrt{2\pi}} \rmd t$ be the standard norm distribution function, then
$$
\sup\limits_{x\in R}  |F(x) -\Phi (x)| \leq 10\sum_{n=1}^\infty \Enum [|X_n|^3].
$$
\end{lem}

\noindent\textbf{Proof of Theorem \ref{thm: CLT}.\;}
First we show $\displaystyle \frac{R^*_t - \mu (t)}{\sigma_{t}} \Dto N (0,1)$ as $t \to \infty$. By  \eqref{eq: dec4R}, we have 
$$
\displaystyle \frac{R^*_t - \mu (t)}{\sigma_{t}}=\sum_{i=1}^\infty \frac{1_{\{N^{(i)}_t \geq 1\}}-\Pnum (N^{(i)}_t \geq 1)}{\sigma_{t}},
$$
which is a sum of independent random variables with $\Var(\frac{R^*_t - \mu (t)}{\sigma_{t}})=1$. Furthermore we have 
$$
\sum_{i=1}^\infty \Enum \Bigl[ |1_{\{N^{(i)}_t \geq 1\}} -\Pnum (N^{(i)}_t \geq 1)|^3\Bigr] \le \sum_{i=1}^\infty \Enum \Bigl[ |1_{\{N^{(i)}_t \geq 1\}} -\Pnum (N^{(i)}_t 
\geq 1)|^2\Bigr]=\Var(R^*_{t})= \sigma_{t}^2.
$$
Let $F_t(x)$ be the cumulative distribution function of $\frac{R^*_t - \mu (t)}{\sigma_{t}}$. Then using Lemma \ref{lem: CLT}, we obtain
$$
\sup\limits_{x\in R} |F_t(x) -\Phi (x)| \leq \frac{10 \sum_{i=1}^\infty \Enum \bigl[ |1_{\{N^{(i)}_t \geq 1\}} -\Pnum (N^{(i)}_t \geq 1)|^3\bigr]  }{\sigma_{t}^3}\le \frac{10}{\sigma_{t}}.
$$
Since $\sigma_{t}\to  \infty$, we have $\lim\limits_{t\to \infty} \sup\limits_x |F_t(x) -\Phi (x)| =0$, proving $\displaystyle \frac{R^*_t - \mu (t)}{\sigma_{t}} \Dto N (0,1)$.

Now Lemma \ref{lem: Difference2CLT} tells us $\displaystyle \frac{R_n - \mu (n)}{\sigma_{n}} \Dto N (0,1)$.
\qed

\section{Proofs of Theorems \ref{thm: MDP} and \ref{thm: LDP}}\label{sec： proof4MDP&LDP}
We have seen that the  proof of   Theorem \ref{thm: CLT} becomes easy  since $R_n^*-R_n$ is negligible in view of Lemma \ref{lem: Difference2CLT}. 
Therefore we want a similar statement which says that there is little difference between $R_n$ and $R_n^*$ so that the results of MDP and LDP for $R^*_n$ can be smoothly transfered into those for $R_n$.

Let's begin with a  well known Chernoff bound for Poisson distribution whose proof is omitted.
  \begin{lem}\label{lem: Poisson-estimates}
Let $\lam>0$ and $X \sim \mathrm{Poisson} (\lam)$. Then
\begin{eqnarray}
\label{eq: Poisson-eq-2}
\Pnum (X \geq x) &\leq& e^{-\lam  } \cdot \bigl(\frac {e\lam}{x}\bigr)^x, \,\,\, \forall\,
x>\lam  ;\\
\label{eq: Poisson-eq-3}
\Pnum (X \leq x) &\leq& e^{-\lam  } \cdot \bigl(\frac {e\lam}{x}\bigr)^x, \,\,\, \forall\,
0\le x<\lam  .
\end{eqnarray}
\end{lem}

By Lemma \ref{lem: Poisson-estimates}, immediately we have the following estimates for $\Pnum (N_{n-w_n} \geq n)+\Pnum (N_{n+w_n} \leq n)$.
\begin{lem}\label{lem: Poisson-2}
Let $h(n)$ and $w_n$ be positive such that $h(n) \to \infty$,  $\frac {w_n}{n} \to 0$  and  $\frac{w^2_n}{nh(n)} \to \infty$. Then
\begin{equation}
\lim_{n \to \infty} \frac{1}{h (n)}\log \bigl[ \Pnum (N_{n-w_n} \geq n)+\Pnum (N_{n+w_n} \leq n) \bigr]=-\infty.
\end{equation}
\end{lem}

 The following lemma  plays an important role in the proofs of Theorems \ref{thm: MDP} and \ref{thm: LDP}.  
\begin{lem}\label{lem: final}
Let $h(n)$ be positive such that $\displaystyle \lim_{n \to \infty} h(n)=\infty$ and $\displaystyle \varliminf_{n \to \infty} \frac {\mu(n)}{h(n)}>0$.
Then 
\begin{equation}\label{eq: MDP-**}
\lim_{n \to \infty} \frac{1}{h (n)}\log\Pnum(|R^*_{n}-R_{n}| \geq \vep  \sqrt{\mu (n) h(n)}) = -\infty, \; \forall \vep>0.
\end{equation}
\end{lem}
\Proof
Fix $\vep>0$. Let $x_n=\vep  \sqrt{\mu (n) h(n)}$ and  $w_n=\Bigl[\frac {n^3h^2(n)}{\mu(n)}\Bigr]^{1/4}$. By 
the definition of $w_n$, as $n\to \infty$
$$
\frac {w_n}{n}=\Bigl[\frac{\mu(n)}{n}\Bigr]^{1/4} \cdot \sqrt{\frac{h(n)}{\mu (n)}} \to 0, \quad \frac{w^2_n}{nh(n)}
=\sqrt{\frac{n}{\mu (n)}} \to \infty.
$$
So, $n\ge 2w_n$ for all sufficiently large $n$. Similar to \eqref{e:minuR*Rscale}, we have
\begin{equation}\label{e:minuR*Rscale10}
\Pnum (|R^*_{n}-R_{n}| \geq x_n) \leq \Pnum (N_{n+w_n} \leq n)+\Pnum (N_{n-w_n} \geq n)+\Pnum (R^*_{n+w_n} -R^*_{n-w_n} 
\geq x_n).
\end{equation}
By Lemma \ref{lem: Poisson-2} we have
\begin{equation}\label{e:N_n_tail_loglim}
\lim\limits_{n\to \infty}\frac{1}{h(n)}\log \bigl[ \Pnum (N_{n+w_n} \leq n)+\Pnum (N_{n-w_n} \geq n) \bigr]
=-\infty.
\end{equation}
For any fixed $\lam  >0$, by the Markov inequality  we can write
\begin{align*} 
&\Pnum (R^*_{n+w_n} -R^*_{n-w_n} \geq x_n)\\
\le& e^{-\lam x_n}  \Enum \Bigl[ e^{\lam ( R^*_{n+w_n} -R^*_{n-w_n}  )} \Bigr] =  e^{-\lam x_n} \Enum \Bigl[ 
\exp\{\lam \sum_{i=1}^\infty 1_{\{N^{(i)}_{n-w_n}=0, N^{(i)}_{n+w_n} -N^{(i)}_{n-w_n}
\geq 1\}} \}\Bigr]\\
=&  e^{-\lam x_n} \prod_{i=1}^\infty \Enum \Bigl[ \exp\{\lam 1_{\{N^{(i)}_{n-w_n}=0, N^{(i)}_{n+w_n} -N^{(i)}_{n-w_n}
\geq 1\}} \} \Bigr].
\end{align*}
Hence
$$
\log \Pnum (R^*_{n+w_n} -R^*_{n-w_n} \geq x_n)\le -\lam   x_n+\sum_{i=1}^\infty \log \Enum \Bigl[ \exp\{\lam 1_{\{N^{(i)}_{n-w_n}=0, 
N^{(i)}_{n+w_n} -N^{(i)}_{n-w_n} \geq 1\}} \} \Bigr].
$$
Since $\Pnum( N^{(i)}_{n-w_n}=0, N^{(i)}_{n+w_n} -N^{(i)}_{n-w_n} \geq 1)=  e^{-\pi_n (n-w_n)}(1- e^{-2\pi_n w_n})$, an easy calculation reveals
$$
\log\Enum \Bigl[ \exp\{\lam \cdot  1_{\{N^{(i)}_{n-w_n}=0, N^{(i)}_{n+w_n} -N^{(i)}_{n-w_n}
\geq 1\}} \} \Bigr] \le (e^\lam-1) \cdot \Bigl[ e^{-\pi_{_i}  (n-w_n)}-e^{-\pi_{_i}  (n+w_n)} \Bigr],
$$
which implies
\begin{align*}
\log \Pnum (R^*_{n+w_n} -R^*_{n-w_n} \geq x_n)\le& -\lam x_n+\sum_{i=1}^\infty (e^\lam-1) \cdot \Bigl[ e^{-\pi_{_i}  (n-w_n)}
-e^{-\pi_{_i}  (n+w_n)} \Bigr]\\
=&-\lam x_n+(e^\lam-1) \cdot \Bigl[ \mu(n+w_n) -\mu(n-w_n) \Bigr].
\end{align*}

Now we choose $\lam=\lam_n :=\ln \bigl( 1+[\frac{n}{\mu(n)}]^{1/4} \bigr)$. The above inequality can be rewritten as 
\begin{equation}\label{eq: prob-bound-*}
\log \Pnum (R^*_{n+w_n} -R^*_{n-w_n} \geq x_n)\le -\lam_n x_n+ \bigl[ \frac{n}{\mu(n)} \bigr]^{1/4} \cdot \Bigl[ 
\mu(n+w_n) -\mu(n-w_n) \Bigr]. 
\end{equation}
Since  $\ddot{\mu}(t)<0, \dot{\mu}(n-w_n)\le \frac{\mu(n-w_n)}{n-w_n}\le \frac{2\mu(n)}{n}$, we have  
$$
\mu(n+w_n)-\mu(n-w_n)\le 2w_n \dot{\mu}(n-w_n)\le \frac{4w_n\mu(n)}{n}.
$$
Substituting $x_n=\vep  \sqrt{\mu (n) h(n)}$ and $w_n=\Bigl[\frac {n^3h^2(n)}{\mu(n)}\Bigr]^{1/4}$ into \eqref{eq: prob-bound-*}, we then get
$$
\frac{1}{h(n)} \log \Pnum (R^*_{n+w_n} -R^*_{n-w_n} \geq x_n)\le - \sqrt{\frac{\mu (n)}{h(n)}} \cdot  \Bigl( \vep \cdot \lam_n -4 \Bigr).
$$
Noting $\lam_n \to \infty$ and $\displaystyle \varliminf_{n \to \infty} \frac {\mu(n)}{h(n)}>0$, we clearly have
\begin{equation}\label{e:mainPart_R_log}
\lim\limits_{n\to \infty}\frac{1}{h(n)}\log \Pnum (R^*_{n+w_n} -R^*_{n-w_n} \geq x_n)=-\infty.
\end{equation}
Combining \eqref{e:minuR*Rscale10}, \eqref{e:N_n_tail_loglim} and \eqref{e:mainPart_R_log}, we  complete the proof of the lemma. 
\qed

In order to prove Theorems \ref{thm: MDP} and \ref{thm: LDP}, we introduce some notations. Set $Y^*_{t}=\frac{R_{t}^*}{\mu(t)}-1$ for $t\ge 1$. Let $\lam \in \Rnum$ and write  
\begin{equation}\label{eq: def4Lambda}
\Lam_t (\lam) := \log \Enum [e^{\lam Y^*_{t}}]
\end{equation} 
for the cumulant generating function associated with $Y^*_{t}$.  Since $\{N_\cdot^{(i)}\}_{i=1}^\infty$ are independent,
\begin{align*}
\Lam_t (\lam) =&\log \Enum \Bigl[ \exp\bigl\{ \frac{\lam}{\mu  (t)}   (R^*_t-\mu  (t)) \bigr\} \Bigr]
= \log \Enum \Bigl[ \exp \bigl\{ \frac{\lam}{\mu (t)} \sum_{i=1}^\infty (e^{-t \pi_{_i} } -1_{\{N^{(i)}_t=0\}}) \bigr\} \Bigr] \\
=&\sum_{i=1}^\infty \log \Enum \Bigl[ \exp \bigl\{\frac{\lam}{\mu  (t)}  (e^{-t \pi_{_i} } -1_{\{N^{(i)}_t=0\}}) \bigr\} \Bigr].
\end{align*}
Since $\Pnum(N_t^{(i)}=0)=e^{-t\pi_{_i} }$, we have
\begin{eqnarray*}
&& \log \Enum \Bigl[ \exp \bigl\{ \frac{\lam}{\mu  (t)} (e^{-t \pi_{_i} } -1_{\{N^{(i)}_t=0\}}) \bigr\} \Bigr] \\
&=&\log \Bigl[ e^{-t \pi_{_i} } \exp \bigl\{ \frac{\lam}{\mu  (t)} (e^{-t \pi_{_i} } -1) \bigr\} +(1-e^{-t \pi_{_i} }) \exp \bigl\{
\frac{\lam}{\mu  (t)} e^{-t \pi_{_i} } \bigr\} \Bigr]\\
&=& \frac{\lam}{\mu  (t)} \cdot (e^{-t \pi_{_i} } -1)+\log \Bigl[ 1+(e^{\frac{\lam}{\mu  (t)}}-1) (1-e^{-t \pi_{_i} }) \Bigr].
\end{eqnarray*}
Combining these equations with $\displaystyle \mu(t)=\sum_{i=1}^\infty (1-e^{-t\pi_{_i} })$, we obtain
\begin{equation}\label{eqno: 3.2}
\Lam_t (\lam) = -\lam +\sum_{i=1}^\infty \log \Bigl[ 1+(e^{\frac{\lam}{\mu  (t)}}-1)   (1-e^{-t \pi_{_i} }) \Bigr].
\end{equation}

\noindent\textbf{Proof of Theorem \ref{thm: MDP}.\;}
Let $\pi$ be regular with index $\gamma\in (0,1)$. Let $\lam  \in R$. Fix  an increasing function $b(t) \to \infty$ 
with $\frac{\mu  (t)}{b (t)} \to +\infty$. Write $u_t= \sqrt{\frac{b (t)}{\mu  (t)}}$ for short. Then $u_t\to 0$ as 
$t \to \infty$. Let 
$$
Z^*_t=\sqrt{\frac{\mu  (t)}{b (t)}} \cdot \Bigl[ \frac{R^*_t}{\mu  (t)} -1 \Bigr]=\frac{R^*_t -\mu  (t)}{\sqrt{\mu  (t) b (t)}},
$$
and let $\wt{\Lam}_t (\lam)=\log\Enum [ e^{\lam   Z_t^*}]$
be the cumulant generating function associated with $Z_t^*$. Noting $Z^*_t= \frac{Y^*_{t}}{u_t}$ and \eqref{eq: def4Lambda} 
we have $\wt{\Lam}_t (\lam  u_t^2 \mu(t))  = \Lam_t (\lam  u_t \mu(t))$. Then by \eqref{eqno: 3.2}, we can write
$$
\frac{\wt{\Lam}_t (\lam  b (t))}{b (t)}=\frac{\wt{\Lam}_t (\lam  u_t^2 \mu(t))}{u_t^2 \mu(t)} = -\frac{\lam}{u_t}
+\frac{1}{u_t^2 \mu(t)} \sum_{i=1}^\infty \log \bigl[ 1+(e^{\lam  u_t}-1)   (1-e^{-t \pi_{_i} }) \bigr].
$$
Note that  $\displaystyle \sum_{i=1}^\infty (1-e^{-t \pi_{_i} }) = \mu  (t)$, $\displaystyle \sum_{i=1}^\infty (1-e^{-t \pi_{_i} })^2
= 2 \mu  (t)-\mu  (2 t)$ and $\displaystyle \sum_{i=1}^\infty (1-e^{-t \pi_{_i} })^3 \leq \mu  (t)$. Noting also $\lim
\limits_{t\to \infty}u_t= 0$ and the fact $\log (1+x)=x-\frac{x^2}{2}+O(x^3)$ as  $x \to 0$, we can write further
\begin{eqnarray*}
\frac{\wt{\Lam}_t (\lam b (t))}{b (t)} &=& -\frac{\lam}{u_t}+ \frac{1}{u_t^2 \mu(t)} \Bigl[ (e^{\lam u_t}-1)\mu(t)
-\frac{1}{2} (e^{\lam u_t}-1)^2 (2\mu(t)-\mu(2t))+O((e^{\lam u_t}-1)^3\mu(t)) \Bigr]\\
&=& -\frac{\lam}{u_t}+ \frac{1}{u_t^2 } \Bigl[ (e^{\lam u_t}-1) -  (e^{\lam u_t}-1)^2 ( 1-\frac{\mu(2t)}{2\mu(t)})
+O((e^{\lam u_t}-1)^3) \Bigr].
\end{eqnarray*}
By $e^{\lam   u_t}=1+\lam u_t+\frac{\lam^2}{2}u_t^2+o(u_t^2)$ as $t\to \infty$, $\frac{\wt{\Lam}_t (\lam b (t))}{b (t)}
=\frac{\lam^2}{2} \cdot \bigl[ \frac{\mu(2t)}{\mu(t)}-1 \bigr]+o(1)$. Since $\pi$ is regular with index $\gamma\in (0,1)$, 
\eqref{e:mu estimate} implies $\mu(t)=\Gamma(1-\gamma) \zeta(t)[1+o(1)]$ and
$$
\lim_{t \to \infty} \frac{\wt{\Lam}_t (\lam b (t))}{b (t)} =\lim_{t \to \infty} \frac{\lam^2}{2} \cdot \bigl[ 
\frac{\zeta(2t)}{\zeta(t)}-1 \bigr] = (2^{^\gamma} -1) \cdot \frac{\lam^2}{2}.
$$
In view of \cite[Theorem 2.3.6]{D&Z},
\begin{equation}\label{thm: MDP_*1-2}
\begin{array}{rcl}
-\frac{1}{2 (2^{^\gamma} -1)} \inf \limits_{x \in \Gamma^o} x^2 \leq& \varliminf \limits_{n \to \infty} \frac{1}{b (n)} \log \Pnum (Z^*_n \in \Gamma)&\\
\leq& \varlimsup \limits_{n \to \infty} \frac{1}{b (n)} \log \Pnum (Z^*_n \in \Gamma)&\leq -\frac{1}{2 (2^{^\gamma} -1)} \inf  \limits_{x \in \bar{\Gamma}} x^2.
\end{array}
\end{equation}

Recall $Z_n=\frac{R_n -\mu (n)}{\sqrt{\mu (n)   b (n)}  }$. By \eqref{eq: MDP-**} and noting $Z^*_n-Z_n=\frac{R_n^* 
-R_n}{\sqrt{\mu (n) b (n)}}$,  
$$
\lim_{n\to \infty} \frac{1}{b(n)} \log \Pnum(|Z_n^*-Z_n|\ge \vep )=-\infty.
$$ 
Hence \eqref{thm: MDP_*1-2} still holds when $Z_n^*$ is replaced  by $Z_n$.
\qed\\

\noindent\textbf{Proof of Theorem \ref{thm: LDP}.\;}
Let $\zeta$ be the $\gamma$-regular function associated with $\pi$ as defined in \eqref{eq:def-phi}. Fix $\lam>0$.  
By \eqref{eqno: 3.2}, we can formulate
$$
\frac{\Lam_t(\lam \mu(t))}{\mu(t)}=-\lam+ \frac{1}{\mu(t)} \sum_{i=1}^\infty \log \bigl[ 1+(e^\lam-1)(1-e^{-t \pi_{_i} }) 
\bigr].
$$
For any fixed $\vep \in (0, 0.1)$, there exists $n_0$ such that $\pi_n\le  \frac{1+\vep }{\zeta^{-1}(n)}, \forall n> n_0$. 
And we have 
\begin{align*}
\frac{\Lam_t(\lam \mu(t))}{\mu(t)}\le & -\lam +\frac{1}{\mu(t)} \Bigl(2n_0 C_\lam  +\sum_{i=1}^\infty 
\log \Bigl[ 1+(e^\lam  -1)(1 -e^{-\frac{(1+\vep )t}{\zeta^{-1}(i)}}) \Bigr] \Bigr),
\end{align*}
where $C_\lam  =\max\limits_{s\in[0,1]} \log \bigl[ 1+(e^\lam  -1)s \bigr]$. In view of  $\mu(t)=\Gamma(1-\gamma)
\zeta(t)[1+o(1)]$ and Lemma \ref{l:Lambda_0.9}, we have 
$$
\lim_{t\to \infty} \frac{1}{\zeta((1+\vep )t)}\sum_{i=1}^\infty \log \bigl[ 1+(1-e^{-\frac{(1+\vep ) t}{\zeta^{-1}(i)}})
(e^{\lam  }-1) \bigr] =\gamma\int_0^\infty \log \bigl[ 1+(1-e^{-s})(e^{\lam  }-1) \bigr]   \frac{\rmd s}{ s^{1+ \gamma}}.
$$
By the regular property of $\zeta$, we have $\lim\limits_{t\to \infty}\frac{\zeta((1+\vep )t)}{\zeta(t)}
=(1+\vep )^\gamma$. So, 
$$
\varlimsup_{t\to \infty}\frac{\Lam _t(\lam   \mu(t))}{\mu(t)} \le  -\lam  +\frac{(1+\vep )^\gamma \gamma }{\Gamma(1-\gamma)}
\int_0^\infty \log \bigl[ 1+(e^\lam  -1)(1-e^{-s}) \bigr] \frac{\rmd s}{s^{1+\gamma}}.
$$

Similarly, by  $\pi_n\ge \frac{1-\vep }{\zeta^{-1}(n)}$ for all sufficiently large $n$, we can prove
$$
\varliminf_{t\to \infty}\frac{\Lam _t(\lam   \mu(t))}{\mu(t)} \ge  -\lam  +\frac{(1-\vep )^\gamma \gamma }{\Gamma 
(1-\gamma)}\int_0^\infty \log \bigl[ 1+(e^\lam  -1)(1-e^{-s}) \bigr] \frac{\rmd s}{s^{1+\gamma}}.
$$
By letting $\vep  \downarrow 0$ we then obtain for each $\lam  > 0$,
\begin{equation}\label{e:LDP_mmmm}
\lim_{t\to \infty}\frac{\Lam_t(\lam \mu(t))}{\mu(t)} =  -\lam +\frac{\gamma}{\Gamma(1-\gamma)}\int_0^\infty 
\log \bigl[ 1+(e^\lam  -1)(1-e^{-s}) \bigr]  \frac{\rmd s}{s^{1+\gamma}}.
\end{equation}

In the same spirit, one can  prove the validity of \eqref{e:LDP_mmmm} for $\lam   \leq 0$. By \eqref{e:LDP_mmmm}
and \eqref{eq: MDP-**}, and exploiting \cite[Theorem 2.3.6]{D&Z}, we finish the proof of  Theorem \ref{thm: LDP} as that of Theorem \ref{thm: MDP}.
\qed\\

\noindent \textbf{Acknowledgement} This work was supported by the National Natural Science Foundation of China (Grant Nos.
 11790273 and  12271351). 


\end{document}